\documentclass[reqno]{amsart}
\usepackage{hyperref}
\usepackage{latexsym}

\newtheorem{theorem}{Theorem}[section]
\newtheorem{lemma}[theorem]{Lemma}
\newcommand{\bt}{\begin{theorem}}
\newcommand{\et}{\end{theorem}}
\newcommand{\blem}{\begin{lemma}}
\newcommand{\elem}{\end{lemma}}

\numberwithin{equation}{section}

\def\open#1{\setbox0=\hbox{$#1$}
\baselineskip = 0pt
\vbox{\hbox{\hspace*{0.4 \wd0}\tiny $\circ$}\hbox{$#1$}} 
\baselineskip = 11pt\!}
\def\fn{\open{f}}
\def\uno{\open{u}_1}
\def\unt{\open{u}_2}
\def\unon{\open{u}_{1,n}}
\def\untn{\open{u}_{2,n}}

\newcommand{\ep}{\hspace*{\fill}$\Box$}
\newcommand{\eps}{\varepsilon}
\newcommand{\pr}{\noindent{\bf Proof. }}

\newcommand{\R}{\mathbb R}
\newcommand{\N}{\mathbb N}

\newcommand{\cs}{\ensuremath{{\mathcal C}} }

\newcommand{\beast}{\begin{eqnarray*}}
\newcommand{\eeast}{\end{eqnarray*}}

\newcommand{\al}{\alpha}

\newcommand{\vphi}{\varphi}

\newcommand{\pa}{\partial}
\newcommand{\ti}{\tilde}

\newcommand{\beq}{ \begin{equation} }
\newcommand{\eeq}{\end{equation} }
\newcommand{\bea}{\begin{eqnarray}}
\newcommand{\eea}{\end{eqnarray}}
\newcommand{\beas}{\begin{eqnarray*}}
\newcommand{\eeas}{\end{eqnarray*}}
\newcommand{\beqs}{\begin{equation*}}
\newcommand{\eeqs}{\end{equation*}}

\newcommand{\ben}{\begin{enumerate}}
\newcommand{\een}{\end{enumerate}}
\newcommand{\ba}{\begin{array}}
\newcommand{\ea}{\end{array}}
\newcommand{\brem}{\begin{thr} {\bf Remark. }\rm}
\newcommand{\ethi}{\end{thr}}

\newcommand{\linfn}[1]{\|#1\|_\infty}
\newcommand{\lon}[1]{\|#1\|_1}
\newcommand{\lt}[1]{\|#1\|_2}
\newcommand{\lonk}[1]{\|#1\|_\mathrm{kin}}
\newcommand{\supp}{\mathop{\mathrm{supp}}}
\renewcommand{\div}{\mathop{\mathrm{div}}}

\begin{document}

\title[Global Weak Solutions of the VKG System]%
{Global Weak Solutions of the Relativistic Vlasov-Klein-Gordon System}

\author[M. Kunzinger]{Michael Kunzinger}
\address{Institut f\"ur Mathematik\\
Strudlhofgasse 4\\ 1090 Wien\\ Austria}
\email{Michael.Kunzinger@univie.ac.at}
\urladdr{http://www.mat.univie.ac.at/\~{}mike/}
\author[{G. Rein}]{Gerhard Rein}
\address{Institut f\"ur Mathematik\\
Strudlhofgasse 4\\ 1090 Wien\\ Austria}
\email{Gerhard.Rein@univie.ac.at}
\urladdr{http://www.mat.univie.ac.at/\~{}rein/}
\author[R. Steinbauer]{Roland Steinbauer}
\address{Institut f\"ur Mathematik\\
Strudlhofgasse 4\\ 1090 Wien\\ Austria}
\email{Roland.Steinbauer@univie.ac.at}
\urladdr{http://www.mat.univie.ac.at/\~{}stein/}
\author[G. Teschl]{Gerald Teschl}
\address{Institut f\"ur Mathematik\\
Strudlhofgasse 4\\ 1090 Wien\\ Austria\\ and International Erwin Schr\"odinger
Institute for Mathematical Physics, Boltzmanngasse 9\\ 1090 Wien\\ Austria}
\email{Gerald.Teschl@univie.ac.at}
\urladdr{http://www.mat.univie.ac.at/\~{}gerald/}
\thanks{Partially supported by  the the Austrian Science Fund's 
Wittgenstein 2000 Award of P.~A.~Markowich.}

\maketitle
                                                                           
\begin{abstract}
We consider an ensemble of classical particles coupled to a
Klein-Gordon field. For the resulting nonlinear system
of partial differential equations, which we call the relativistic
Vlasov-Klein-Gordon system, we prove the existence of global weak
solutions for initial data satisfying a size restriction.
The latter becomes necessary since the energy of the system
is indefinite, and only for restricted data a-priori bounds
on the solutions can be derived from conservation of energy.
\end{abstract}

\section{Introduction}\label{intro}
\setcounter{equation}{0}

When considering the interaction of classical particles with
classical or quantum fields various different situations arise:
On the one hand one can consider the coupling of a single classical particle
to a field. How this should properly be done for the case of a Maxwell
field is a classical problem, cf.\ \cite{Ab},
and the effective dynamics and asymptotics
of such systems is an active field of research, 
cf.\ \cite{IKS,KKS1,KKS2,KS} and the references
there. On the other hand, in kinetic theory ensembles of classical particles
are considered which interact by fields which they create collectively. 
There is an extensive literature on such systems, with particles interacting
by non-relativistic, gravitational or electrostatic 
fields --- the Vlasov-Poisson system ---, by electrodynamic 
fields --- the Vlasov-Maxwell system ---, or by general relativistic 
gravity --- the Vlasov-Einstein system. In all these systems the only
interaction of the particles is via the fields which they create
collectively, a situation which is sometimes referred to as the mean 
field limit of a many-particle system.

In the present paper we consider an ensemble of particles which can move
at relativistic speeds and interact by a quantum mechanical Klein-Gordon
field. Let $f=f(t,x,v)\geq 0$ denote the density of the particles
in phase space, $\rho=\rho(t,x)$ their density in space, and
$u=u(t,x)$ a scalar Klein-Gordon field; $t \in \R$, $x \in \R^3$, and
$v \in \R^3$ denote time, position, and momentum respectively.
The system then reads as follows: 
\begin{equation} \label{vl}
\pa_t f + \hat v \cdot \pa_x f - \pa_x u \cdot \pa_v f = 0,
\end{equation}
\begin{equation} \label{kg}
\pa_t^2 u - \Delta u + u = - \rho,
\end{equation}
\begin{equation} \label{rhodef}
\rho(t,x) = \int f(t,x,v)\,dv.
\end{equation}
Here we have set all physical constants as well as the rest mass
of the particles to unity, and
\begin{equation} \label{hatvdef}
\hat v = \frac{v}{\sqrt{1+|v|^2}}
\end{equation}
denotes the relativistic velocity of a particle with momentum $v$.
We call this system the relativistic Vlasov-Klein-Gordon system.
 
To our knowledge it has not yet been considered in the literature. 
Our motivation for initiating a study of this system
is the following. In \cite{IKM} a single classical particle coupled to
a Klein-Gordon field is considered, and the system 
(\ref{vl}), (\ref{kg}), (\ref{rhodef}) is meant as a natural 
generalization of this to the many-particle situation.
On the other hand, the system falls within the general class of
nonlinear PDE systems like the Vlasov-Poisson or Vlasov-Maxwell system,
and by studying it one may hope to learn more about the general
properties of this important class of problems from mathematical physics.
A major issue in this area, which we also focus on in the present paper,
is the question of global existence of solutions to the corresponding
initial value problem. 

The coupling in the system above is set up in such a way that the system
is conservative;
\begin{equation} \label{energy}
\int \sqrt{1 + |v|^2}f\,dx\,dv + 
\frac{1}{2} \int \bigl[|\pa_t u|^2 + 
|\pa_x  u|^2 + |u|^2\bigr]\,dx + \int \rho u \,dx =:
E_K + E_F + E_C
\end{equation}
is conserved along sufficiently regular solutions.
From experience with the related systems from kinetic theory mentioned above
one knows that for global existence questions it is essential to
derive a-priori bounds on the solutions from conservation of energy.
In this context the relativistic Vlasov-Klein-Gordon system
poses the following specific difficulty: The interaction term
$E_C$ in the energy need not be positive. Therefore, one has to try to estimate
this term in terms of the positive quantities $E_K$ and $E_F$.
By H\"older's inequality, interpolation, and Sobolev's inequality,
\[
\left| \int \rho u \,dx \right|
\leq
\|\rho(t)\|_{6/5} \|u(t)\|_6 \leq C E_K(t)^{1/2} \|\pa_x u(t)\|_2; 
\]
the details of a more general version of this estimate can be found below,
cf.\ Lemma~\ref{enest}. The problem is that, after applying the Cauchy
inequality, the right hand side in this estimate is of the same order
of magnitude as the positive terms in the energy, and no a-priori bound
on the solution seems to follow from conservation of energy. 
This situation is similar to the gravitational case of the relativistic
Vlasov-Poisson system where it is known that spherically symmetric
solutions with negative energy blow up in finite time \cite{GlSch1}.
Our way out of this difficulty is to observe that in the estimate above
the constant $C$ depends on $\|\fn\,\|_1$ and $\|\fn\,\|_\infty$,
$\fn$ being the initial datum for $f$,
and if we require that it is smaller than an appropriate threshold then
a-priori bounds on the solutions can be derived from conservation of energy.
These bounds are such that one can pass to the limit along a sequence
of global solutions to an appropriately regularized system,
and this limit is a weak solution to the relativistic Vlasov-Klein-Gordon
system. The details of these arguments are then similar to the corresponding 
ones for the Vlasov-Maxwell system \cite{DL}, cf.\ also \cite{KR}. 

The paper proceeds as follows: In the next section we collect for easier
reference some results on the linear, inhomogeneous Klein-Gordon equation.
In Section~3 we prove a global existence and uniqueness result
for smooth solutions to an appropriately regularized version of
the Vlasov-Klein-Gordon system. We then proceed to derive uniform a-priori
bounds on a sequence of such regularized solutions from conservation
of energy for restricted data, and prove the existence of a global, 
weak solution to the original system. This is done in Section~4, where we 
also discuss some of the properties of these weak solutions.

We conclude this introduction with some remarks on related results
in the literature; we restrict ourselves to results on the initial
value problem in the three dimensional case. 
For the Vlasov-Poisson system global classical
solutions for general initial data have been established in
\cite{Pf,LP,Sch}. For the relativistic Vlasov-Maxwell system
such global classical solutions are so far only known for
small, nearly neutral, or nearly spherically symmetric data
\cite{GlStr,GlSch2,R1}. For general data the existence
of global weak solutions was obtained in \cite{DL}.
When passing from the Vlasov-Poisson to the Vlasov-Maxwell
system a lot of the difficulties with classical solutions
of course arise from the different field equation which becomes
hyperbolic instead of elliptic. However, it must be emphasized
that already for the relativistic Vlasov-Poisson system
where the only difference to the Vlasov-Poisson system is that a $v$ in the
Vlasov equation is replaced by a $\hat v$ no global existence
result is known for classical solutions with general initial data,
and for the gravitational case where the energy is indefinite
classical solutions can blow up in finite time as mentioned above,
cf.\ \cite{GlSch1}. An investigation of the fully relativistic
gravitational situation, i.e., of the Vlasov-Einstein system,
was initiated in \cite{RR1}.   

\section{The linear Klein-Gordon equation}
\setcounter{equation}{0}

Although our notation is mostly standard or self-explaining
we explicitly mention the following conventions:
For a function $h=h(t,x,v)$ or $h=h(t,x)$ we denote for given $t$
by $h(t)$ the corresponding function of the remaining variables.
By $\|\,.\,\|_p$ we denote the usual $L^p$-norm for $p\in[1,\infty]$. 
The index $c$ in function spaces refers to compactly supported
functions. Throughout the paper
the convolution denoted by $\ast$ refers to the spatial variables. 

For easier reference we now collect some results on the 
linear, inhomogeneous Klein-Gordon equation
\begin{equation} \label{linkg}
\pa_t^2 u  - \Delta u  + u = g 
\end{equation}
with a prescribed right hand side $g=g(t,x)$. For $t>0$ and $x\in \R^3$ let
\[
R(t,x) := \frac{1}{4 \pi} \frac{\delta(|x|-t)}{t} - \frac{1}{4 \pi}
H(t-|x|) \frac{J_1(\sqrt{t^2-|x|^2})}{\sqrt{t^2-|x|^2}}
\]
where $\delta$ is the $\delta$-distribution, $H$ is the Heaviside function
and $J_1$ is the Bessel function of the first kind.

Given initial data
\begin{equation} \label{indalinkg}
u(0) = \uno,\
\pa_t u(0) = \unt
\end{equation}
the solution of the initial value problem (\ref{linkg}), (\ref{indalinkg})
can be written as
\[
u(t,x) = u_{\mathrm{hom}}(t,x) + u_{\mathrm{inh}}(t,x),\ t\geq 0,\ x\in \R^3.
\]
Here
\begin{eqnarray*}
u_{\mathrm{hom}}(t,x)
&:=&
(\pa_t R \ast \uno)(t,x) + (R\ast \unt)(t,x) \\
&=&
\frac{1}{4 \pi t^2} \int_{|y|=t} \uno(x-y)\, dS_y
- \frac{1}{4 \pi t^2} \int_{|y|=t} \nabla \uno(x-y)\cdot y\, dS_y\\ 
&&
{}
- \frac{1}{8 \pi} \int_{|y|=t} \uno(x-y)\, dS_y
- \frac{1}{4 \pi} \int_{|y|\leq t} \uno(x-y)\, 
\left(\frac{J_1(\xi)}{\xi}\right)'
\frac{t}{\xi} dy\\ 
&&
{}
+ \frac{1}{4 \pi t} \int_{|y|=t} \unt(x-y)\, dS_y
- \frac{1}{4 \pi} \int_{|y|\leq t} \unt(x-y)\, 
\frac{J_1(\xi)}{\xi} dy
\end{eqnarray*}
with the abbreviation
\[
\xi:= \sqrt{t^2 - |y|^2}
\]
is the solution of the homogeneous Klein-Gordon equation (\ref{linkg})
with initial data (\ref{indalinkg}),
and
\begin{eqnarray*}
u_{\mathrm{inh}}(t,x)
&:=&
\int_0^t \int R(t-s,x-y)\, g(s,y)\,dy\,ds \\
&=& 
\frac{1}{4 \pi}\int_0^t \int_{|x-y|=t-s}\!\!\!\!\!\!
g(s,y) \,dS_y\frac{ds}{t-s} 
- \frac{1}{4 \pi}\int_0^t \int_{|x-y|\le t-s}\!\!\!\!\!\! 
g(s,y)\,\frac{J_1(\xi)}{\xi} dy\,ds
\end{eqnarray*}
with the abbreviation
\[
\xi:=\sqrt{(t-s)^2 - |x-y|^2}
\]
is the solution of the inhomogeneous Klein-Gordon equation with vanishing
initial data. These formulas can be found in \cite{MS} and \cite{Sid}.
They can be established by observing that the substitution
$w(t,x,\zeta) := u(t,x)\exp(-i\zeta)$ transforms (\ref{linkg}) into 
a wave equation for $w$ which can be solved in the usual way.

If $\al \in \N_0^{\,3}$ denotes an arbitrary multi-index and $\pa^\al$
the corresponding spatial derivative the above formulas imply the 
following estimate for the solution of (\ref{linkg}), (\ref{indalinkg}),
which is not very sophisticated but good enough for our purpose:
\begin{equation}\label{kgestimate}
\|\pa^\al u(t)\|_\infty \le C (1+t)^4
\left(\linfn{\pa^\al \uno} + \linfn{\nabla \pa^\al  \uno}
+ \linfn{\pa^\al \unt}
+ \linfn{\pa^\al g(t)}\right).
\end{equation}

\section{Global classical solutions of the regularized system} \label{regsys}
\setcounter{equation}{0}

The aim of the present section is to establish the following global 
existence and 
uniqueness result for a suitably regularized Vlasov-Klein-Gordon system:
\bt \label{rvkgdeltath} 
Let $\delta \in \mathcal{C}_c^\infty(\R^3)$. Consider the regularized,
relativistic Vlasov-Klein-Gordon system where the right hand side 
$-\rho$ in (\ref{kg}) is replaced by $- \rho* \delta$.
For initial data $0 \le \fn \in \mathcal{C}_c^1(\R^6)$, 
$\uno\in \mathcal{C}^3_b(\R^3)$ and
$\unt\in \mathcal{C}^2_b(\R^3)$
there exists
a unique solution $(f,u)$ to the regularized system 
with $f\in \mathcal{C}^1([0,\infty[\times \R^6)$ 
and $u\in \mathcal{C}^2([0,\infty[\times \R^3)$,
satisfying the initial conditions $f(0)=\fn$, $u(0)=\uno$, $\pa_t u(0)=\unt$.
\et
\pr 
We begin by recursively defining iterates 
$f_n: [0,\infty[ \times \R^6 \to [0,\infty[$,
$u_n: [0,\infty[ \times \R^3 \to \R$:
Let  $f_0(t,z):=\fn(z)$, $z=(x,v)$, 
and suppose $f_n$ has already been defined.
Let $\rho_n(t,x):= \int f_n(t,x,v)\,dv$ and let $u_n$ 
denote the solution to the Klein-Gordon
equation (\ref{linkg}) with right hand side $- \rho_n * \delta$ and 
initial data $\uno$, $\unt$, cf.\ the previous section.
Denote by $Z_n(s,t,z)=(X_n,V_n)(s,t,x,v)$ the solution 
of the characteristic system 
\beast
\dot x &=& \hat v, \\
\dot v &=& - \pa_x u_n(s,x)
\eeast
with initial datum $Z_n(t,t,z) = z$. The next iterate is then defined as 
$f_{n+1}(t,z):= \fn(Z_n(0,t,z))$, i.e., as the solution of the Vlasov
equation (\ref{vl}) with $u$ replaced by $u_n$ and initial datum $\fn$.

The flow of the characteristic system
is measure preserving, hence
$\lon{f_n(t)} = \lon{\rho_n(t)} = \lon{\fn\,}$. This implies that
for any $\al \in \N_0^{\,3}$, $\linfn{\pa^\al (\rho_n(t)\ast \delta)}$
is bounded, uniformly in $t\geq 0$ and $n\in\N$.
Let $T>0$ be arbitrary. In what follows constants denoted by $C$
may change from line to line and depend on the initial data, on 
the regularization kernel $\delta$, and on $T$, but never on $n\in \N$.
By (\ref{kgestimate}), $\linfn{\pa^\alpha u_n(t)}\leq C$ for
$|\al|\leq 2$ and $t\in [0,T]$. Hence there exist $R>0$ and $P>0$ such that 
$f_n(t,x,v) = 0$ if $|v| > P$ or $|x| > R$ and $t\in [0,T]$, 
in particular, $f_n(t)\in \mathcal{C}^1_c(\R^6)$.  
Using these bounds it is now easy to show that $(f_n)_{n\in\N}$
is a uniform Cauchy sequence on $[0,T]\times \R^6$: 
By definition,
\[
\linfn{f_{n+1}(t) - f_n(t)} \le 
C \linfn{Z_n(0,t,\,.\,) - Z_{n-1}(0,t,\,.\,)}.
\]
Using the characteristic system and the bound on $\pa_x^2 u_n(t)$ we obtain, 
abbreviating $Z_n(s,t,x,v)$ by $Z_n(s)$,
\begin{eqnarray*}
|X_n(s) - X_{n-1}(s)| 
&\le& 
\int_s^t |V_n(\tau) - V_{n-1}(\tau)|\,d\tau  \\
|V_n(s) - V_{n-1}(s)| 
&\le& 
C \int_s^t |X_n(\tau) - X_{n-1}(\tau)|\,d\tau \\
&&
{}+
\int_s^t \linfn{\pa_x u_n(\tau) -\pa_x u_{n-1}(\tau)}\,d\tau, 
\end{eqnarray*}
hence
\[
|Z_n(0) - Z_{n-1}(0)| \le C
\int_0^t |Z_n(\tau) - Z_{n-1}(\tau)|\,d\tau
+ \int_0^t \linfn{\pa_x u_n(\tau) -\pa_x u_{n-1}(\tau)}\,d\tau.
\]
By Gronwall's inequality, 
\[
\linfn{f_{n+1}(t) - f_n(t)} 
\le C \int_0^t \linfn{\pa_x u_n(\tau) -\pa_x u_{n-1}(\tau)}\,d\tau .
\]
On the other hand, the formulas for $u_n$ from the previous section,
the fact that we have regularized the right hand side of (\ref{kg}),
and the above bounds on the support of $f_n$ imply that
\begin{eqnarray*} 
\linfn{\pa_x u_n(\tau) -\pa_x u_{n-1}(\tau)} 
&\leq& 
C \lon{\rho_n(\tau) - \rho_{n-1}(\tau)} \\
&\leq& 
C \int\!\!\!\!\int |f_n(\tau,x,v) - f_{n-1}(\tau,x,v)|\,dv\,dx \\
&\leq&
C \linfn{f_n(\tau) - f_{n-1}(\tau)}.
\end{eqnarray*}
Combining these estimates we obtain 
\[
\linfn{f_{n+1}(t) - f_{n}(t)} \le C \int_0^t 
\linfn{f_n(\tau) - f_{n-1}(\tau)}\,d\tau.
\]
Thus, by induction,
\[
\linfn{f_{n+1}(t) - f_{n}(t)} \le C \frac{C^n t^n}{n!}\,,
\]
so $(f_n)_{n\in \N}$ is uniformly Cauchy on $[0,T] \times \R^6$. 
The same is true for $(\rho_n)_{n\in \N}$ 
and $(\pa^\al u_n)_{n\in\N}$ for $|\al|\leq 2$.
That the uniform limit $(f,u)$ of the iterative
sequence has the required regularity and is the unique solution
of our initial value problem on $[0,T]$
follows, and since $T>0$ was arbitrary the proof is complete.
\ep

In the next section we want to obtain a global weak solution of the
relativistic Vlasov-Klein-Gordon system as a limit of a sequence
of solutions to systems regularized with $\delta_n$'s that converge to
the $\delta$-distribution. To do so we will need energy bounds on
these regularized solutions, but since we have modified the system
energy conservation takes a somewhat different form from what was stated
in the introduction:
\blem \label{conservation}
Let $d\in \cs^\infty_c(\R^3)$ be even, $\delta = d*d$,
and $\uno$, $\unt \in \cs_c(\R^3)$. 
Let $0\le \fn \in \mathcal{C}_c^1(\R^6)$ and 
let $(f,u)$ be the unique solution to the regularized system according to
Theorem~\ref{rvkgdeltath} 
with initial conditions $f(0) = \fn$,
$u(0) = \uno*\delta$, $\pa_t u(0) = \unt*\delta$. 
Let $\rho(t,x) = \int f(t,x,v)\, dv$
and denote by $\tilde u$ the unique solution to the initial value problem
\[
\pa_t^2 \tilde u - \Delta \tilde u + \tilde u = -\rho* d\,,
\]
\[
\tilde u(0) = \uno*d,\ 
\pa_t \tilde u(0) = \unt*d . 
\]
Then
\[
\tilde E := \int \sqrt{1+|v|^2}f\,dx\,dv + 
\frac{1}{2} \int [|\pa_t \tilde u|^2 + 
|\pa_x \tilde u|^2 + |\tilde u|^2]\,dx + \int \rho u \,dx =: 
E_K + \tilde E_F + E_C
\]
is constant in $t$.
\elem
\pr 
Using the Vlasov equation and integration by parts we obtain
\[
\frac{d}{dt} E_K = -\int \pa_t \rho \, u\, dx\,.
\]
Also, 
\[
\frac{d}{dt} \tilde E_F = -\int \pa_t \tilde u \, \rho* d\, dx = 
-\int \pa_t(\tilde u*d)\,\rho\,dx;
\]
for the last equality observe that $d$ is assumed to be even.
The convolution
$\tilde u *d$ satisfies the Klein-Gordon equation with right hand side 
$-\rho*\delta$ and 
initial conditions $(\tilde u*d)(0) = \uno*\delta$, 
$\pa_t(\tilde u*d)(0) = \unt*\delta$. 
Hence by uniqueness, $u=\tilde u*d$. Summing up,
\[
\frac{d}{dt} E_C = \int \pa_t u\, \rho\,dx + \int u\,\pa_t\rho \,dx =
- \frac{d}{dt} \tilde E_F - \frac{d}{dt} E_K.
\]
\ep

\section{Weak solutions}
\setcounter{equation}{0}

Based on Theorem \ref{rvkgdeltath} we now prove existence of 
global weak solutions to 
the relativistic Vlasov-Klein-Gordon system. 
The following auxiliary result will allow us to derive a-priori bounds
from conservation of energy, at least for appropriately restricted
initial data:
\blem \label{enest} 
Let $p \in ]3/2,\infty]$ and $1/p + 1/q = 1$.
In addition to the assumptions of Lemma~\ref{conservation}
we assume further that $d\geq 0$ with $\int d =1$.
Let $(f,u)$ be a solution as obtained in 
Lemma \ref{conservation}. Then 
\[
\left| \int u(t,x)\, \rho(t,x)\, dx \right|
\leq C(\fn\,)\, \|\pa_x \tilde u(t)\|_2\, E_K(t)^{1/2},\ t\geq 0,
\]
where
\[
C(\fn\,) :=  \left(\frac{4q}{\pi}\right)^{1/2} 3^{-7/6} 
\left(\frac{q+3}{q}\right)^{(q+3)/6}\,
\lon{\fn\,}^{(3-q)/6} \|\fn\,\|_p^{q/6}.
\]
\elem
\pr 
For any $R>0$,
\[
\rho(t,x) = \int_{|v|\le R} f \,dv + \int_{|v|>R} f\,dv \le 
\left(\frac{4\pi}{3}R^3\right)^{1/q} \|f(t,x,\,.\,)\|_p + 
\frac{1}{R} \int \sqrt{1+|v|^2} f\,dv.
\]
If we choose $R$ such that the right hand side becomes minimal we obtain
the estimate 
\[
\rho(t,x) \leq C_q
\|f(t,x,\,.\,)\|_p^{q/(q+3)} 
\left(\int \sqrt{1+|v|^2} f \, dv\right)^{3/(q+3)}
\]
where
\[
C_q:=\left(\frac{4 \pi}{3}\right)^{1/(q+3)}
\frac{q+3}{3} \left(\frac{3}{q}\right)^{q/(q+3)}.
\]
We take this estimate
to the power $(q+3)/(q+2)$, integrate in $x$, apply
H\"older's inequality, and observe that $\|f(t)\|_p=\|\fn\,\|_p$ to obtain
\[
\|\rho(t)\|_{(q+3)/(q+2)} \le C_q \|\fn\,\|_p^{q/(q+3)} E_K(t)^{3/(q+3)}.  
\]
By assumption on $p$ we have $(q+3)/(q+2) > 6/5$, and by interpolation,
\begin{equation} \label{rho6/5}
\|\rho(t)\|_{6/5} \le \|\rho(t)\|_1^{(3-q)/6}
\|\rho(t)\|_{(q+3)/(q+2)}^{(q+3)/6}
\leq C_q^{(q+3)/6} \lon{\fn\,}^{(3-q)/6}
\|\fn\,\|_p^{q/6} E_K(t)^{1/2}.
\end{equation}
By H\"older's inequality, Young's inequality, and Sobolev's inequality,
\beast
\left|\int \rho u\,dx\right| 
&=& 
\left|\int \rho*d\, \tilde u\, dx \right|
\le \|\tilde u(t)\|_6 \|\rho(t)*d\|_{6/5} \le \frac{1}{\sqrt{S_3}}
\lt{\pa_x\tilde u(t)}  \|\rho(t)\|_{6/5}\\
&\leq&
\frac{1}{\sqrt{S_3}} C_q^{(q+3)/6} \lon{\fn\,}^{(3-q)/6}
\|\fn\,\|_p^{q/6} \lt{\pa_x\tilde u(t)} E_K(t)^{1/2}
\eeast
with $S_3 := 3 (\pi/2)^{4/3}$, cf.~\cite[8.3]{LL};
recall from the proof of Lemma~\ref{conservation} that $u=\tilde u*d$.
\ep

We now proceed to the main result of the present paper. 
In its formulation we employ the space
\[
L^1_\mathrm{kin}(\R^6):= 
\left\{f: \R^6 \to \R \mid \|f\|_\mathrm{kin}:= 
\int\!\!\!\!\int \sqrt{1+|v|^2} |f(x,v)|\,dx\,dv < \infty\right\}.
\]
\bt \label{main} 
Let $\fn\in L^1_\mathrm{kin}(\R^6)\cap L^p(\R^6)$ for some $p\in[2,\infty]$,
$\fn \ge 0$, 
$\uno \in H^1(\R^3)$, $\unt \in L^2(\R^3)$, and assume that
\[
\lon{\fn\,}^{(3-q)/3} \|\fn\,\|_p^{q/3} <
\frac{\pi}{2q} 3^{7/3} \left(\frac{q}{q+3}\right)^{(q+3)/3}
\]
where $1/p+1/q=1$.
Then there exists a global weak
solution $(f,u)$ of the relativistic Vlasov-Klein-Gordon system
(\ref{vl}), (\ref{kg}), (\ref{rhodef}) with these initial data, 
more precisely, 
\[
f\in L^\infty([0,\infty[\times L^p(\R^6)),\
u\in L^\infty([0,\infty[\times H^1(\R^3))
\]
with 
\[
\pa_t u\in L^\infty([0,\infty[\times L^2(\R^3))
\]
such that the following holds:
\begin{itemize}
\item[{\rm (a)}]
$(f,u)$ satisfies (\ref{vl}), (\ref{kg}), (\ref{rhodef}) in 
$\mathcal{D}'(]0,\infty[\times \R^6)$.
\item[{\rm (b)}]
The mapping $[0,\infty[\ni t \mapsto (f(t),u(t),\pa_t u(t)) 
\in L^2(\R^6)\times L^2(\R^3)\times L^2(\R^3)$
is weakly continuous with $(f,u,\pa_t u)(0)=(\fn,\uno,\unt)$.
\item[{\rm (c)}]
$f(t) \geq 0$ a.e., $\|f(t)\|_p \leq \|\fn\|_p,\ t\geq 0$,
and $\pa_t \rho + \div j = 0$ in $\mathcal{D}'(]0,\infty[\times \R^3)$
where $j(t,x):= \int \hat v f(t,x,v)\, dv$.
The weak solution conserves mass: 
$\|f(t)\|_1 = \|\fn\|_1$ for a.~a.\ $t\geq 0$. 
\end{itemize}
\et

\noindent
{\bf Remark.} Regardless of whether one considers data satisfying the
restriction of the theorem as small or not it should be emphasized
that this is not a small data result of the type known for example
for the Vlasov-Maxwell system, cf.\ \cite{GlStr}, since it does 
not rely on the fields being small and corresponding dispersive
effects of the free streaming Vlasov equation. In particular,
it is worthwhile to note that there is no size restriction
on the data for the Klein-Gordon field.    

\noindent
{\bf Proof of Theorem~\ref{main}.}
By our assumption on $\fn$, $C(\fn\,)^2 < 2$ where the left hand 
side is defined as in Lemma~\ref{enest}. We choose $\eps \in ]0,1[$
such that $C(\fn\,)^2 < 2\eps$.
Next we choose sequences $(\fn_n)$ in $\mathcal{C}^\infty_c(\R^6)$,
$\fn_n \ge 0$, $(\unon)$, $(\untn)$ in 
$\mathcal{C}^\infty_c(\R^3)$ 
such that
$\fn_n\to \fn$ in $L^1_\mathrm{kin}(\R^6) \cap L^p(\R^6)$, 
$\unon \to \uno$
in $H^1(\R^3)$, $\untn \to \unt$ in $L^2(\R^3)$,
and we require that
\begin{equation}\label{indarestrn}
\sup_{n\in \N} C(\fn_n)^2 < 2\eps.
\end{equation}
Let $d_n\in \mathcal{C}_c^\infty(\R^3)$ be
non-negative, even, of unit integral, and with 
$\supp d_n \subseteq B_{1/n}(0)$, and define $\delta_n := d_n*d_n$.
Denote by $(f_n, u_n)$ the regularized solution according to
Theorem~\ref{rvkgdeltath} with $\delta_n$ 
replacing $\delta$ and initial
conditions
\[
f_n (0) = \fn_n,\ 
u_n(0) = \unon*\delta_n,\
\pa_t u_n(0) = \untn*\delta_n.
\]
Let $\rho_n := \int f_n\, dv$ and denote by $\tilde u_n$ the solution 
to the Klein-Gordon equation 
with right hand side $ - \rho_n* d_n$ and initial data
\[
\tilde u_n(0) = \unon* d_n,\
\pa_t \tilde u_n(0) = \untn* d_n.
\]
As in the proof of Lemma~\ref{conservation} it follows that 
$u_n = \tilde u_n * d_n$. 
Using Lemma~\ref{enest} and Cauchy's 
inequality we find that
\beast
\tilde E 
&\geq&
\|f_n(t)\|_\mathrm{kin} + \frac{1}{2}\lt{\pa_t \tilde u_n(t)}^2 + 
\frac{1}{2}\lt{\tilde u_n(t)}^2 + \frac{1}{2}\lt{\pa_x \tilde u_n(t)}^2 
\\
&&{}- C(\fn_n) 
\lt{\pa_x \tilde u_n(t)}\,\|f_n(t)\|_\mathrm{kin}^{1/2}\\
&\geq&
(1-\eps)\|f_n(t)\|_\mathrm{kin} + \frac{1}{2}\lt{\pa_t \tilde u_n(t)}^2 + 
\frac{1}{2}\lt{\tilde u_n(t)}^2
+ \left(\frac{1}{2} - \frac{1}{4\eps} \sup_{k\in \N}C(\fn_k)^2\right)
\lt{\pa_x \tilde u_n(t)}^2.
\eeast 
Hence by conservation of energy and (\ref{indarestrn})
there exists $C>0$ such that 
for all $n \in \N$ and $t\geq 0$,
$\lonk{f_n(t)}$, $\lt{\pa_t u_n(t)}$, $\|u_n(t)\|_{H^1}\leq C$. Also, 
$\lt{f_n(t)}=\lt{\fn_n}$ is bounded by interpolation. 
Since the bounds are independent
of $n$ and $t$, by extracting appropriate subsequences 
(again denoted by the same indices)
we obtain $f\in L^\infty([0,\infty[;L^2(\R^6))$, 
$u\in L^\infty([0,\infty[; H^1(\R^3))$ with 
$\pa_t u\in L^\infty([0,\infty[; L^2(\R^3))$ and 
\begin{eqnarray*}
f_n 
&\rightharpoonup& 
f \mbox{ in } L^2([0,T]\times \R^6),\\
u_n 
&\rightharpoonup& 
u \mbox{ in } L^2([0,T]\times \R^3),\\
\pa_x u_n 
&\rightharpoonup& 
\pa_x u \mbox{ in } L^2([0,T]\times \R^3),\\
\pa_t u_n 
&\rightharpoonup& 
\pa_t u \mbox{ in } L^2([0,T]\times \R^3)
\end{eqnarray*}
for all $T>0$.
Moreover, $f\in L^\infty([0,\infty[;L^1_\mathrm{kin}(\R^6))$ since
\[
\int \sqrt{1+v^2} f \,dx\,dv = 
\lim_{R\to\infty}\lim_{n\to\infty} 
\int_{{|x|\le R} \atop |v|\le R} \sqrt{1+v^2} f_n \,dx\,dv 
\le C\,.
\]
Since $\|f_n(t)\|_p = \|\fn_n\|_p$ is bounded as well, 
(\ref{rho6/5}) implies boundedness
of $\|\rho_n(t)\|_{6/5}$. Hence without loss of generality we may assume 
that $\rho_n \rightharpoonup
\rho$ in $L^{6/5}([0,T]\times \R^3)$ for any $T$. Moreover, $\|u_n(t)\|_6 
\le C \|\pa_x u_n(t)\|_2$, so we may finally
suppose that $u_n \rightharpoonup u$ in $L^6([0,T]\times \R^3)$ for any $T$.
We claim that $u$ satisfies the Klein-Gordon equation in the weak sense. 
Indeed, let $\vphi\in \mathcal{C}^\infty_c(]0,\infty[\times \R^3)$. Then
\begin{eqnarray*}
0 
&=& 
\int \left(\pa_t^2 u_n - \Delta u_n + u_n + \rho_n*\delta_n\right)
\vphi\,dx\,dt \\
&=& 
\int \left(u_n \pa_t^2 \vphi - u_n \Delta \vphi + u_n\vphi + 
(\rho_n*\delta_n)\vphi\right) \,dx\,dt\\
&\to& 
\int \left(u\pa_t^2 \vphi - u \Delta \vphi + u\vphi + \rho\vphi\right) \,dx\,dt
\end{eqnarray*}
which yields the claim; note that 
$\int (\rho_n\ast\delta_n)\vphi = \int \rho_n (\delta_n\ast \vphi)$
and $\delta_n\ast \vphi \to \vphi$ in $L^6(\R^3)$.

Turning now to the Vlasov equation, we adapt an argument from 
\cite[Sec.~3]{KR}. 
For $\eps>0$ let $\xi_\eps \in \mathcal{C}_c^\infty(\R)$, $0\le \xi_\eps
\le 1$, $\xi_\eps = 1$ on $[\eps,T]$, $\supp \xi_\eps \subseteq [\eps/2,2T]$. 
Then
\[
\pa_t(f_n\xi_\eps) + \hat v \pa_x(f_n\xi_\eps) = 
\mathrm{div}_v(f_n \xi_\eps \pa_x u_n) + f_n\xi_\eps'
\]
Hence by the velocity-averaging lemma of Golse, Lions, Perthame and Sentis 
(\cite{GLPS, DL}) 
we have 
\begin{eqnarray*}
&\forall R> 0 \ \forall \psi\in \mathcal{C}_c^\infty(B_R(0)) \ 
\exists C=C(R,\psi) \ \forall n\in \N,\, 
\eps> 0 \,: &\\
&\left| \int \xi_\eps(\,.\,) f_n(\,.\,,\,.\,,v) \psi(v)\,dv \right| 
\in H^{1/4}(\R\times\R^3)& \\
&\| \int \xi_\eps(\,.\,) f_n(\,.\,,\,.\,,v) \psi(v)\,dv \|_{H^{1/4}} 
\le C \left(\lt{\xi_\eps f_n}^2
+ \lt{f_n \xi_\eps \pa_x u_n}^2  + \lt{f_n\xi_\eps'}^2\right)^{1/2}\,.&
\end{eqnarray*}
In the case $p=\infty$ the boundedness properties derived above already 
assure the boundednes of the sequence $(\int \xi_\eps(\,.\,) f_n(\,.\,,\,.\,,v) 
\psi(v)\,dv)_{n\in \N}$ in $H^{1/4}(\R\times \R^3)$. The case $p<\infty$ is
more delicate and will be treated below.
Since the restriction operator from $H^{1/4}(\R\times\R^3)$ to
$L^2([0,T]\times B_{R}(0))$ is compact, by a diagonal sequence argument, 
for each $\psi\in \mathcal{C}_c^\infty(\R^3)$ we may extract
a subsequence (again denoted by $(f_n)_{n\in \N}$) independent of 
$\eps:= 1/m$, $R'=m$ ($m\in \N$)
such that
\begin{equation}\label{l2conv}
\int  f_n(\,.\,,\,.\,,v) \psi(v)\,dv \to 
\int  f(\,.\,,\,.\,,v) \psi(v)\,dv
\end{equation}
in $L^2(]0,T[\times B_{R'}(0))$.

For $p<\infty$ we show the validity of (\ref{l2conv})
as follows: For $\eta>0$ we define 
$\beta_\eta(\tau):=\tau/(1+\eta\tau)$. Then $\beta_\eta\circ f_n$,
in addition to satisfying Vlasov's equation (with force $\partial_xu_n$)
and the same bounds as $f_n$, is bounded in $L^\infty$, and we may use the
previous case. In the limit $\eta\to 0$ we obtain $L^1$-convergence in  
(\ref{l2conv}) which together with the uniform integrability of 
$\{(\int f^n(.,.,v)\psi(v)dv)^2\}$ implies the desired $L^2$-convergence.

We use this property to prove that $f$ is a weak solution 
of the Vlasov equation. For
$\vphi_1\in \mathcal{C}_c^\infty(\R)$, $\vphi_2,\, \vphi_3 
\in \mathcal{C}_c^\infty(\R^3)$,
\begin{eqnarray*}
0
&=& 
\int_0^T \int\!\!\!\!\int
f_n(t,x,v) \Bigl(\vphi_1'(t)\vphi_2(x)\vphi_3(v) +
\hat v\cdot \pa_x \vphi_2(x) \vphi_1(t)\vphi_3(v) \\
&&
\phantom{\int_0^T \int\!\!\!\!\int\Bigl(\ }
- \pa_x u_n(t,x)\cdot \pa_v\vphi_3(v)\vphi_1(t)\vphi_2(x)\Bigr)\,dv\,dx\,dt
\end{eqnarray*}
Choose a subsequence as above for $\psi = \pa_v \vphi_2$ and $R'>0$ 
such that $\supp \vphi_2 \subseteq
B_{R'}(0)$. Then 
\[
\int f_n(\,.\,,\,.\,,v)\pa_v\vphi_3(v)\,dv \to \int f(\,.\,,\,.\,v)\pa_v\vphi_3(v)\,dv
\]
in $L^2([0,T]\times B_{R'}(0))$ and $\pa_x u_n
\vphi_1\vphi_2
\rightharpoonup \pa_x u\vphi_1\vphi_2$ in 
$L^2((0,T)\times B_{R'}(0))$.
Thus, finally,
\begin{eqnarray*}
0
&=& 
\int_0^T \int\!\!\!\!\int f(t,x,v) \Bigl(\vphi_1'(t)\vphi_2(x)\vphi_3(v) 
+ \hat v\cdot \pa_x \vphi_2(x) \vphi_1(t)\vphi_3(v) \\
&& 
\phantom{\int_0^T \int\!\!\!\!\int f(t,x,v) \Bigl(\ }
 - \pa_x u(t,x)\cdot \pa_v\vphi_3(v)\vphi_1(t)\vphi_2(x)\Bigr)\,dv\,dx\,dt
\end{eqnarray*}
which establishes (a).

As to (b), we integrate the Vlasov equation once with respect to $t$ and define
\[
\tilde{f}(t):= \fn - \int_0^t \hat{v} \cdot\pa_x f(s) ds - 
\int_0^t \pa_x u(s) \cdot \pa_v f(s) ds
\]
(to be understood as an identity between distributions in $x$ and $v$).
Then, using the construction of $f$ it is not hard to see that $f(t)=\ti{f}(t)$ in $\mathcal{D}'(\R^6)$
for a.e.\ $t\in[0,\infty[$, that is, we can replace $f$ by $\ti{f}$.
Moreover, $\ti{f}(0)=\fn$ and a straightforward application of the dominated
convergence theorem shows that $t\mapsto \ti{f}(t)$ is weakly continuous.
A similar argument applies to $u$.

As to (c), we only prove that mass is conserved, the other assertions being
standard. Let $\eps >0$ and choose
$R>0$ such that $\int_{|x|\geq R} \rho (0,x)\,dx < \eps$.
Without loss of generality we can assume the same for $\rho_n(0)$,
but since for the regularized problem the particles move along 
characteristics with speeds $|\hat v| < 1$ we have
$\int_{|x|\geq R+T} \rho_n (t,x)\,dx < \eps$ on any time interval $[0,T]$.
Hence,
\beast
\int \rho (t,x)\, dx
&\geq&
\int_{|x|\leq R+T} \rho (t,x)\,dx 
= \lim_{n\to \infty} \int_{|x|\leq R+T} \rho_n (t,x)\,dx\\
&\geq&
\lim_{n\to \infty} \int \rho_n (t,x)\,dx - \eps
= 
\|\fn\|_1 - \eps.
\eeast
Since on the other hand $\int \rho(t,x)\, dx \leq \|\fn\|_1$ and
$\eps$ is arbitrary this proves conservation of mass. 
\ep
\medskip\\
\noindent {\bf Acknowledgments.} 
We thank Peter A.\ Markowich for several helpful discussions.

\end{document}